\documentclass[12pt,reqno]{amsart}
\usepackage{enumerate, latexsym, amsmath, amsfonts, amssymb, amsthm, graphicx, color}
\usepackage{hyperref}
\hypersetup{hypertex=true,
	colorlinks=true,
	linkcolor=blue,
	anchorcolor=blue,
	citecolor=blue}
\usepackage{amsmath}
\def\pmod #1{\ ({\rm{mod}}\ #1)}
\def\i{\mathrm{i}}

\def\Q{\Bbb Q}

\def\l{\left}
\def\r{\right}
\def\bg{\bigg}
\def\({\bg(}
\def\){\bg)}
\def\t{\text}
\def\f{\frac}

\def\adj{{\rm adj}}

\def\ls{\le}

\def\ve{\varepsilon}

\def\eq{\equiv}

\def\u{{\bf u}}
\def\v{{\bf v}}
\def\Proof{\noindent{\it Proof}}

\theoremstyle{plain}
\newtheorem{theorem}{Theorem}

\newtheorem{lemma}{Lemma}

\theoremstyle{definition}

\theoremstyle{remark}
\newtheorem{remark}{Remark}

\newcommand{\Legendre}[2]{\left(\frac{#1}{#2}\right)}
\makeatletter
\@namedef{subjclassname@2010}{%
	\textup{2010} Mathematics Subject Classification}
\makeatother
\vspace{4mm}

\allowdisplaybreaks[1]
\begin{document}	
	\hbox{Preprint}
\medskip

	\title[A determinant involving Legendre symbols]
{On a determinant involving linear combinations of Legendre symbols}
	
	\author[K. Liu, Z.-W. Sun and L.-Y. Wang]{Keqin Liu*, Zhi-Wei Sun and Li-Yuan Wang}
	
	\address {(Keqin Liu) School of Mathematics and Physics, Xi'an Jiaotong-Liverpool University,
	Suzhou, 215123, People's Republic of China} \email{Keqin.Liu@xjtlu.edu.cn}

\address {(Zhi-Wei Sun) School of Mathematics, Nanjing University\\
	Nanjing, 210093, People's Republic of China} \email{zwsun@nju.edu.cn}

\address {(Li-Yuan Wang) School of Physical and Mathematical Sciences, Nanjing Tech University, Nanjing, 211816, People's Republic of China}
	\email{\tt lywang@njtech.edu.cn}

	\begin{abstract}
   In this paper, we prove a conjecture of the second author by evaluating the determinant
  $$\det\left[x+\left(\frac{i-j}p\right)+\left(\frac ip\right)y+\left(\frac jp\right)z+\left(\frac{ij}p\right)w\right]_{0\ls i,j\ls(p-3)/2}$$
  for any odd prime $p$, where $(\frac{\cdot}p)$ denotes the Legendre symbol.
  In particular, the determinant is equal to $x$ when $p\equiv 3\pmod4$. 
	\end{abstract}
	
	\thanks{
		\newline \indent
		The three authors Liu, Sun and Wang, were supported by the Leadership Talent Program (Science and Education) of SIP (Grant No. KJQ2024202), the National Natural Science Foundation of China (Grant Nos. 12371004), and the National Natural Science Foundation of China (Grant No. 12201291), respectively.
		\newline\indent{*Corresponding author.}
        }
	\maketitle
	\section{Introduction}	
	\setcounter{lemma}{0}
	\setcounter{theorem}{0}
	\setcounter{corollary}{0}
	\setcounter{remark}{0}
	\setcounter{equation}{0}
	\setcounter{conjecture}{0}

For a matrix $A=[a_{ij}]_{1\ls i,j\ls n}$ over a field, we use $\det A$ or $|a_{ij}|_{1\ls i,j\ls n}$
to denote its determinant. If $a_{ji}=-a_{ij}$ for all $i,j=1,\ldots,n$, then we say that $A$ is skew-symmetric, and in this case we have
$$\det A=|a_{ji}|_{1\ls i,j\ls n}=|-a_{ij}|_{1\ls i,j\ls n}=(-1)^n\det A.$$
In particular, $\det A=0$ if $A$ is skew-symmetric and $n$ is odd.

Let $p$ be an odd prime, and let $(\f{\cdot}p)$ be the Legendre symbol.
For any integer $m$, Lehmer \cite{L} and Carlitz \cite{carlitz} found all the eigenvalues of the matrices
$$\l[x+\l(\f{j+k+m}p\r)\r]_{1\ls j,k\ls p-1}\ \l[x+\l(\f{j-k+m}p\r)\r]_{1\ls j,k\ls p-1}$$
respectively.

Let $\ve_p$ and $h_p$ denote the fundamental unit and class number
	of the real quadratic field $\mathbb{Q}(\sqrt{p})$, respectively.
For convenience, we write
\begin{align}\label{abp}
	 \varepsilon_p^{h_p}=a_p+b_p\sqrt{p}\ \ \t{with}\ a_p,b_p\in\mathbb{Q},
\end{align}
and
\begin{align}\label{ab'p}
	 \varepsilon_p^{(2-(\frac{2}{p}))h_p}=a_p'+b_p'\sqrt{p}\ \ \t{with}\ a_p',b_p'\in\mathbb{Q}.
\end{align}
In 2004, motivated by coding theory, Chapman \cite{chapman} proved that
 \begin{align*}&\ \l|x+\l(\f{i+j-1}p\r)\r|_{1\ls i,j\ls(p+1)/2}
 \\=&\ \begin{cases}(\f 2p)2^{(p-1)/2}(pb_px-a_p)&\t{if}\ p\eq1\pmod4,
 \\2^{(p-1)/2}&\t{if}\ p>3\ \t{and}\ p\eq3\pmod4.\end{cases}
 \end{align*}
 As
 $$\f{p+1}2-i+\f{p+1}2-j-1=p-i-j\eq-(i+j)\pmod p,$$
 we see that
 $$\l|\l(\f{-1}p\r)x+\l(\f{i+j}p\r)\r|_{0\ls i,j\ls(p-1)/2}=\l|x+\l(\f{i+j-1}p\r)\r|_{1\ls i,j\ls(p+1)/2}.$$
 Chapman's ``evil determinant conjecture" (cf. \cite{problems,evil}) states that
	\begin{align}\label{evil}
		 \left|\left(\frac{i-j}{p}\right)\right|_{0\le i,j\le (p-1)/2}=\begin{cases}
		 	-a_p' & \mbox{if}\ p\equiv 1\pmod4,\\
		 	1    & \mbox{if}\ p\equiv 3\pmod4,
		 \end{cases}
	\end{align}
 this was  confirmed by Vsemirnov \cite{V3,V1}  via matrix decomposition.

Sun \cite{S19} defined $M_p$ as the matrix obtained from $[(\f{i-j}p)]_{0\ls i,j\ls (p-1)/2}$ via replacing all the entries in the first row by $1$, and conjectured that
\begin{equation}\label{Mp}\det M_p=\begin{cases}(-1)^{(p-1)/4}&\text{if}\ p\equiv 1\pmod4,
\\(-1)^{(h(-p)-1)/2}&\text{if}\ p>3\ \t{and}\ p\equiv3\pmod4,\end{cases}
\end{equation}
where $h(-p)$ denotes the class number of the imaginary quadratic field $\mathbb Q(\sqrt{-p})$
which is known to be odd. This conjecture was recently confirmed by Wang and Wu \cite{WW}.

Recently, Sun \cite{S24} determined
$$\l|x+\l(\f{i\pm j}p\r)+\l(\f ip\r)y+\l(\f jp\r)z\r|_{0\ls i,j\ls(p-1)/2},$$
and posed many conjectures on determinants involving linear combinations of Legendre symbols.

In this paper, we establish the following result conjectured by Sun \cite[Conjecture 3.2(ii)]{S24}.

	\begin{theorem}\label{Th1} Let $p$ be an odd prime, and let $a_p'$ and $b_p'$ be
given by  $(\ref{ab'p})$.

{\rm (i)} In the case $p\eq1\pmod 4$, we have
\begin{equation}\label{main1}
 \begin{aligned}&\ \l|x+\l(\f{i-j}p\r)+\l(\f ip\r)y+\l(\f jp\r)z+\l(\f{ij}p\r)w\r|_{0\ls i,j\ls(p-3)/2}
\\&\qquad\ =((y+1)(z+1)-wx)\l(\f 2p\r)b_p'-a_p'x.
\end{aligned}
\end{equation}

{\rm (ii)} Provided that $p\eq3\pmod4$, we have
\begin{equation}\label{main3} \l|x+\l(\f{i-j}p\r)+\l(\f ip\r)y+\l(\f jp\r)z+\l(\f{ij}p\r)w\r|_{0\ls i,j\ls(p-3)/2}=x.
\end{equation}
\end{theorem}

It is easy to verify \eqref{main3} for $p=3$. From now on, we let $p$ be a fixed prime greater than $3$.

In the next section, we will reduce the two parts of Theorem \ref{Th1} to
 the identity
\begin{equation}\label{simple}\l|x+\l(\f{i-j}p\r)\r|_{1\ls i,j\ls(p-1)/2}=\begin{cases}
(\frac{2}{p})b'_p-a_p'x&\t{if}\ p\eq1\pmod4,
\\x&\t{if}\ p\eq3\pmod4.
\end{cases}
\end{equation}
This identity in the case $p\eq3\pmod4$ first appeared as a conjecture
in \cite{S19}. In Section 3, we will prove \eqref{simple} in the case $p\eq3\pmod4$.
Based on some lemmas in Section 4, we are going to prove \eqref{simple} for the case $p\eq1\pmod4$
in Section 5.

Throughout this paper, for a matrix $A$ we use $A^T$ to denote the transpose of $A$.
For a matrix $A=[a_{ij}]_{1\ls i,j\ls n}$ over a field, its adjugate matrix is given by $\adj(A)=[A_{ji}]_{1\ls i,j\ls n}$, where $A_{ji}$ is the cofactor of the entry $a_{ji}$ in $A$.

	\section{Reduction of Theorem \ref{Th1} to the identity \eqref{simple}}
	\setcounter{lemma}{0}
	\setcounter{theorem}{0}
	\setcounter{corollary}{0}
	\setcounter{remark}{0}
	\setcounter{equation}{0}
	\setcounter{conjecture}{0}

The following basic lemma can be found in \cite[Lemma 2.1]{RJ}.

\begin{lemma} \label{Lem2.1} Let $A=[a_{ij}]_{0\ls i,j\ls m}$ be a matrix over a field. Then
\begin{equation}\label{ab}\det[x+a_{ij}]_{0\ls i,j\ls m}-\det[a_{ij}]_{0\ls i,j\ls m}
=x\det[b_{ij}]_{1\ls i,j\ls m},
\end{equation}
where $b_{ij}=a_{ij}-a_{i0}-a_{0j}+a_{00}$.
\end{lemma}

Fix an odd prime $p$. Applying Lemma \ref{Lem2.1}, we obtain
\begin{equation}\label{linear}\begin{aligned}&\ \l|x+\l(\f{i-j}p\r)\r|_{1\ls i,j\ls(p-1)/2}=\l|x+\l(\f{i-j}p\r)\r|_{0\ls i,j\ls(p-3)/2}
\\=&\ \l|\l(\f{i-j}p\r)\r|_{0\ls i,j\ls(p-3)/2}
+x\l|\l(\f{i-j}p\r)-\l(\f ip\r)-\l(\f{-j}p\r)\r|_{1\ls i,j\ls(p-3)/2}.
\end{aligned}\end{equation}
On the other hand, by \cite[Theorem 1.2(ii)]{S24} we have
\begin{align*}&\ \l(\f{-1}p\r)\l|x+\l(\f{i-j}p\r)+\l(\f ip\r)y+\l(\f jp\r)z\r|_{0\ls i,j\ls(p-3)/2}
\\=&\ (y+1)\l(z+\l(\f{-1}p\r)\r)\l|\l(\f{i-j}p\r)\r|_{0\ls i,j\ls(p-3)/2}
\\&\ +\l(\f{-1}p\r)x\l|\l(\f{i-j}p\r)-\l(\f ip\r)-\l(\f{-j}p\r)\r|_{1\ls i,j\ls(p-3)/2}.
\end{align*}
Thus, from \eqref{simple} we can deduce the more general conclusion:
\begin{equation}\label{-w}
\begin{aligned}&\ \l|x+\l(\f{i-j}p\r)+\l(\f ip\r)y+\l(\f jp\r)z\r|_{0\ls i,j\ls(p-3)/2}
\\=&\ \begin{cases}(y+1)(z+1)(\f 2p)b_p'-a_p'x&\t{if}\ p\eq1\pmod4,
\\x&\t{if}\ p\eq3\pmod4.
\end{cases}
\end{aligned}
\end{equation}
Combining this with \cite[Theorem 1.2(i)]{S24}, we see that in the case $p\eq1\pmod4$
the equality \eqref{simple} implies that
\begin{equation}\label{no-w}
\begin{aligned}&\ \l|\l(\f{i-j}p\r)+\l(\f ip\r)y+\l(\f jp\r)z+\l(\f{ij}p\r)w\r|_{0\ls i,j\ls(p-3)/2}
\\=&\ \l|\l(\f{i-j}p\r)+\l(\f ip\r)y+\l(\f jp\r)z\r|_{0\ls i,j\ls(p-3)/2}
\\=&\ (y+1)(z+1)\l(\f 2p\r)b_p'.
\end{aligned}
\end{equation}

Now we consider the case $p\eq1\pmod4$. Define
$$A(x,y,z,w)=\l[x+\l(\f{i-j}p\r)+\l(\f ip\r)y+\l(\f jp\r)z+\l(\f{ij}p\r)w\r]_{0\ls i,j\ls(p-3)/2}.$$
If \eqref{simple} holds, then
\begin{gather*}|A(0,0,0,0)|=\l(\f 2p\r)b_p'\not=0,\ |A(1,0,0,0)|=\l(\f 2p\r)b_p'-a_p',
\\ |A(0,1,0,0)|=|A(0,0,1,0)|=2\l(\f 2p\r)b_p'
\ \t{and}\ |A(0,0,0,1)|=\l(\f 2p\r)b_p'
\end{gather*}
with the aids of \eqref{-w} and \eqref{no-w}, hence by applying \cite[Theorem 3.1]{ChenSun}
we obtain
$$|A(x,y,z,w)|=((y+1)(z+1)-wx)\l(\f 2p\r)b_p'-a_p'x.$$

For $i,j=0,\ldots,(p-3)/2$, let
$$a_{ij}=\l(\f{i-j}p\r)+\l(\f ip\r)y+\l(\f jp\r)z+\l(\f{ij}p\r)w.$$
By \cite[Theorem 2.1(i)]{S24}, $\l|a_{ij}\r|_{0\ls i,j\ls (p-3)/2}$
does not depend on $w$.
Combining this with Lemma \ref{Lem2.1}, we get
\begin{align*} &\ \l|x+a_{ij}\r|_{0\ls i,j\ls (p-3)/2}
\\ =&\ \l|\l(\f{i-j}p\r)+\l(\f ip\r)y+\l(\f jp\r)z\r|_{0\ls i,j\ls (p-3)/2}+x|b_{ij}|_{1\ls i,j\ls(p-3)/2},
\end{align*}
where
\begin{align*}b_{ij}&=a_{ij}-a_{i0}-a_{0j}+a_{00}
\\&=a_{ij}-\l(\f ip\r)(y+1)-\l(\f jp\r)\l(z+\l(\f{-1}p\r)\r)+0
\\&=\l(\f{i-j}p\r)-\l(\f ip\r)-\l(\f{-j}p\r)+\l(\f{ij}p\r)w.
\end{align*}
In light of Lemma \ref{Lem2.1},
\begin{align*}|b_{ij}|_{1\ls i,j\ls(p-3)/2}&=\l|\l(\f{ij(i-j)}p\r)-\l(\f jp\r)-\l(\f{-i}p\r)+w\r|_{1\ls i,j\ls(p-3)/2}
\\&=\l|\l(\f{ij(i-j)}p\r)-\l(\f jp\r)-\l(\f{-i}p\r)\r|_{1\ls i,j\ls(p-3)/2}
\\&\quad +w|c_{ij}|_{2\ls i,j\ls(p-3)/2},
\end{align*}
where
\begin{align*}c_{ij}=&\ \l(\f{ij(i-j)}p\r)-\l(\f jp\r)-\l(\f{-i}p\r)
\\&\ -\l(\f{j(1-j)}p\r)+\l(\f jp\r)+\l(\f{-1}p\r)
\\&\ -\l(\f{i(i-1)}p\r)+1+\l(\f{-i}p\r)-\l(\f1p\r)-\l(\f{-1}p\r)
\\=&\ \l(\f{ij(i-j)}p\r)-\l(\f{i(i-1)}p\r)-\l(\f{j(1-j)}p\r).
\end{align*}

Now assume that $p\eq3\pmod4$. Then $|c_{ij}|_{2\ls i,j\ls(p-3)/2}$ vanishes since $[c_{ij}]_{2\ls i,j\ls(p-3)/2}$ is a skew-symmetric matrix of odd order.
Hence $|b_{ij}|_{1\ls i,j\ls(p-3)/2}$ does not depend on $w$.
Therefore
$|x+a_{ij}|_{0\ls i,j\ls(p-3)/2}$ does not depend on $w$, and its value is $x$ provided  \eqref{-w}.

In view of the above, we have reduced Theorem \ref{Th1} to the identity \eqref{simple}.

	\section{Proof of \eqref{simple} in the case $p\eq3\pmod4$}
	\setcounter{lemma}{0}
	\setcounter{theorem}{0}
	\setcounter{corollary}{0}
	\setcounter{remark}{0}
	\setcounter{equation}{0}
	\setcounter{conjecture}{0}

Let $p>3$ be a prime with $p\eq3\pmod4$.
As $[(\f{i-j}p)]_{0\ls i,j\ls(p-3)/2}$
is a skew-symmetric matrix of odd order, we have
\begin{equation}\label{zero}\l|\l(\f{i-j}p\r)\r|_{1\ls i,j\ls(p-1)/2}=\l|\l(\f{i-j}p\r)\r|_{0\ls i,j\ls(p-3)/2}=0.
\end{equation} Combining this with \eqref{linear}, we see that
$$\l|x+\l(\f{i-j}p\r)\r|_{1\ls i,j\ls (p-1)/2}=mx$$
for an integer $m$ not depending on $x$. Thus
$$\l|x+\l(\f{i-j}p\r)\r|_{1\ls i,j\ls (p-1)/2}=x$$
if
\begin{equation}
\label{1+}\l|1+\l(\f{j-i}p\r)\r|_{1\ls i,j\ls (p-1)/2}=1.
\end{equation}
So it suffices to prove \eqref{1+}.

Let $M_p=[m_{ij}]_{0\ls i,j\ls (p-1)/2}$ be the matrix obtained from $[(\f{i-j}p)]_{0\ls i,j\ls(p-1)/2}$
 via replacing all the entries in the first row by $1$, and let $M_p^*=[m_{ij}^*]_{0\ls i,j\ls(p-1)/2}$ be the matrix $[m_{(p-1)/2-i,(p-1)/2-j}]_{0\ls i,j\ls (p-1)/2}$. Then
$$\det M_p^*=\det M_p=(-1)^{(h(-p)-1)/2}$$
by \eqref{Mp}.
For each $j=0,\ldots,(p-1)/2$, we have $$m_{ij}^*=\begin{cases}(\f{j-i}p)&\t{if}\ 0\ls i<(p-1)/2,\\1&\t{if}\ i=(p-1)/2.\end{cases}$$

As $|(\f{j-i}p)|_{0\ls i,j\ls (p-1)/2}=1\not=0$, by Cramer's rule there are rational numbers
$c_0,\ldots,c_{(p-1)/2}$ such that
$$\sum_{i=0}^{(p-1)/2}\l(\f{j-i}p\r)c_i=1\quad\t{for all}\ j=0,\ldots,\f{p-1}2.$$
Define a lower triangular matrix $A=[a_{ij}]_{0\ls i,j\ls(p-1)/2}$ by
$$a_{ij}=\begin{cases}1&\t{if}\ i=j\in\{0,\ldots,(p-3)/2\},
\\c_j&\t{if}\ i=(p-1)/2\ \t{and}\ 0\ls j\ls(p-1)/2,
\\0&\t{otherwise}.\end{cases}$$
Then
$$A\l[\l(\f{j-i}p\r)\r]_{0\ls i,j\ls(p-1)/2}=M_p^*$$
and hence
\begin{equation}
\label{A}\det A=\det A\cdot\l|\l(\f{j-i}p\r)\r|_{0\ls i,j\ls(p-1)/2}=\det M_p^*=(-1)^{(h(-p)-1)/2}.
\end{equation}

Observe that  $M_p^*A^T$ coincides with the matrix $B=[b_{ij}]_{0\ls i,j\ls(p-1)/2}$ with
$$b_{ij}=\begin{cases}(\f{j-i}p)&\t{if}\ i,j\in\{0,\ldots,(p-3)/2\},
\\-1&\t{if}\ i<j=(p-1)/2,
\\1&\t{if}\ i=(p-1)/2>j,\\c&\t{if}\ i=j=(p-1)/2,
\end{cases}$$
where $c=\sum_{i=0}^{(p-1)/2}c_i$. Therefore
\begin{equation}\label{B}\det B=\det M_p^* \cdot \det A=1
\end{equation}
by \eqref{A}.
Let $B_*$ be the matrix obtained from $B$ by replacing the right-bottom entry $c$ by $1$.
Then $\det B_*=\det B=1$ since $|(\f{j-i}p)|_{0\ls i,j\ls(p-3)/2}=0$ by \eqref{zero}.
Via adding the last row of $B_*$ to all previous rows, we see that
$$\det B_*=\l|1+\l(\f{j-i}p\r)\r|_{0\ls i,j\ls(p-3)/2}=\l|1+\l(\f{j-i}p\r)\r|_{1\ls i,j\ls(p-1)/2}.$$
Combining this with \eqref{B}, we immediately get the desired identity \eqref{1+}.

In view of the above, we have completed our proof of \eqref{simple} in the case $p\eq3\pmod4$.

	\section{Some lemmas}
	\setcounter{lemma}{0}
	\setcounter{theorem}{0}
	\setcounter{corollary}{0}
	\setcounter{remark}{0}
	\setcounter{equation}{0}
	\setcounter{conjecture}{0}

To prove \eqref{simple} for primes $p\eq1\pmod4$, we need
	several lemmas.

	\begin{lemma}[The Matrix-Determinant Lemma]\label{mdl}
		Let $H$ be an $m\times m$ matrix over a field $F$, and let $\u,\v$ be two $m$-dimensional column vectors with entries in $F$. Then we have
		\begin{equation*}
			\det (H+\u\v^T)=\det H+ \v^T \adj(H)\u.
		\end{equation*}		
	\end{lemma}
	
\begin{remark} This lemma is well-known. One may consult \cite{mdl} for a  proof.
\end{remark}

The following known lemma can be found in \cite[Theorem 3]{V1}.

\begin{lemma} \label{uv} We have
	\begin{align*}
		&\	\det \bigg[\frac{x_i+y_j}{1+x_iy_j}\bigg]_{1\le i,j\le m}  \\
		= &\ \frac{1}{2}\(\prod_{i=1}^{m}(  1+x_i )(  1+y_i )+(-1)^m \prod_{i=1}^{m}(  1-x_i )(  1-y_i )\)  \\
		&\ \times\prod_{1\le i<j\le m}^{}(x_i-x_j)(y_j-y_i)\cdot\prod_{i=1}^{m}  \prod_{j=1}^{m}(1+x_iy_j)^{-1}.
	\end{align*}
\end{lemma}

\begin{lemma}\label{a+-b} Let $p$ be a prime with $p\eq1\pmod4$, and set $n=\frac{p-1}{2}$ and $\zeta =e^{2\pi \i/p}$.  Then
	\begin{equation}
	  \prod_{j=1}^{n}\( 1+ \(\frac{j}{p}\)\zeta^{-j} \)^2 =(-1)^{n/2}\zeta^{-n(n+1)/2}(b_p' p+ a_p' \sqrt{p}),
\end{equation}
and
\begin{equation}
 \prod_{j=1}^{n}\( 1- \(\frac{j}{p}\)\zeta^{-j} \)^2 =(-1)^{n/2}\zeta^{-n(n+1)/2}(b_p' p- a_p' \sqrt{p}).
	\end{equation}
\end{lemma}
 \begin{proof}
By~\t{\cite[(4.7) and (4.8)]{V1}}, we have
\[
\prod_{j=1}^{n}
\left(
1+\Legendre{j}{p}\zeta^j
\right)^2
=
\zeta^{n(n+1)/2}\varepsilon_p^{\left(
2-\Legendre{2}{p}
\right)h_p}
\prod_{j=1}^{n}
\left(
e^{\pi \i j/p}-e^{-\pi \i j/p}
\right),
\]
and
\[
\prod_{j=1}^{n}
\left(
1-\Legendre{j}{p}\zeta^j
\right)^2
=
\zeta^{n(n+1)/2}\varepsilon_p^{-\left(
2-\Legendre{2}{p}
\right)h_p}
\prod_{j=1}^{n}
\left(
e^{\pi \i j/p}-e^{-\pi \i j/p}
\right).
\]
Moreover, Vsemirnov~\t{\cite[Corollary~2]{V1}} gives
\[
\prod_{j=1}^{n}
\left(
e^{\pi \i j/p}-e^{-\pi \i j/p}
\right)
=
(-1)^{n/2}\sqrt p.
\]
Hence
\begin{align}\label{plus}
\prod_{j=1}^{n}
\left(
1+\Legendre{j}{p}\zeta^j
\right)^2
&=
\zeta^{n(n+1)/2}\varepsilon_p^{\left(
2-\Legendre{2}{p}
\right)h_p}(-1)^{n/2}\sqrt p,
\end{align}
and
\begin{align}\label{minus}
\prod_{j=1}^{n}
\left(
1-\Legendre{j}{p}\zeta^j
\right)^2
&=
\zeta^{n(n+1)/2}\varepsilon_p^{-\left(
2-\Legendre{2}{p}
\right)h_p}(-1)^{n/2}\sqrt p.
\end{align}
By definition,
\[
\varepsilon_p^{\left(
2-\Legendre{2}{p}
\right)h_p}
=
a'_p+b'_p\sqrt p.
\]
It is known that the norm of $\varepsilon_p$ with respect to $\mathbb{Q}(\sqrt{p})/\mathbb{Q}$ is $-1$ and the class number $h_p$ is odd (cf. \t{\cite[pp. 185-187]{Cohn}} or \t{\cite[Lemma~2]{V1}}). 
Note that
\[
\left(
2-\Legendre{2}{p}
\right)h_p
\]
is also odd. Hence
\[
\left(a'_p+b'_p\sqrt p\right)
\left(a'_p-b'_p\sqrt p\right)
=
-1,
\]
so that
\[
\varepsilon_p^{-\left(
2-\Legendre{2}{p}
\right)h_p}
=
-a'_p+b'_p\sqrt p.
\]
Consequently,
\[
\varepsilon_p^{\left(
2-\Legendre{2}{p}
\right)h_p}
+
\varepsilon_p^{-\left(
2-\Legendre{2}{p}
\right)h_p}
=
2b'_p\sqrt p,
\]
and
\[
\varepsilon_p^{\left(
2-\Legendre{2}{p}
\right)h_p}
-
\varepsilon_p^{-\left(
2-\Legendre{2}{p}
\right)h_p}
=
2a'_p.
\]
Adding the two product identities \eqref{plus} and \eqref{minus} therefore gives
\begin{align}\label{bp}
\frac12
\Bigg(
&
\prod_{j=1}^{n}
\left(
1+\Legendre{j}{p}\zeta^j
\right)^2
+
\prod_{j=1}^{n}
\left(
1-\Legendre{j}{p}\zeta^j
\right)^2
\Bigg)
&=
(-1)^{n/2}
\zeta^{n(n+1)/2}
b'_p p.
\end{align}
Subtracting \eqref{minus} from \eqref{plus} gives
\begin{align}\label{ap}
\frac12
\Bigg(
&
\prod_{j=1}^{n}
\left(
1+\Legendre{j}{p}\zeta^j
\right)^2
-
\prod_{j=1}^{n}
\left(
1-\Legendre{j}{p}\zeta^j
\right)^2
\Bigg)
&=
(-1)^{n/2}
\zeta^{n(n+1)/2}
a'_p\sqrt p.
\end{align}

We now pass from $\zeta^j$ to $\zeta^{-j}$. Namely,
\[
1+\Legendre{j}{p}\zeta^{-j}
=
\Legendre{j}{p}\zeta^{-j}
\left(
1+\Legendre{j}{p}\zeta^j
\right),
\]
and
\[
1-\Legendre{j}{p}\zeta^{-j}
=
-\Legendre{j}{p}\zeta^{-j}
\left(
1-\Legendre{j}{p}\zeta^j
\right).
\]
After squaring and multiplying over $1\le j\le n$, both identities
contribute the same factor
\[
\prod_{j=1}^{n}\zeta^{-2j}
=
\zeta^{-n(n+1)}.
\]
Thus
\[
\prod_{j=1}^{n}
\left(
1\pm\Legendre{j}{p}\zeta^{-j}
\right)^2
=
\zeta^{-n(n+1)}
\prod_{j=1}^{n}
\left(
1\pm\Legendre{j}{p}\zeta^j
\right)^2.
\]
Combining this identity with \eqref{bp} and \eqref{ap} gives the two
identities stated in Lemma~4.3.
\end{proof}

	For convenience, we introduce Vsemirnov's notations~\t{\cite{V1}} which will be used soon.

Let $p$ be a prime with $p\eq1\pmod4$, and set $n=\frac{p-1}{2}$ and $\zeta =e^{2\pi i/p}$. Define
\[
C
=
\left[
\left(\frac{j-i}{p}\right)\right]_{0\le i,j\le n}
\]
As in \cite{V1}, we define matrices $D,U, V$ of order $n+1$ whose
$(i,j)$ entries $(0\ls i,j\ls n)$ are as follows:
$$	d_{ij}=\begin{cases}\prod_{0\le k\le n \atop k\ne i}^{}\frac{1}{\zeta^{2i}-\zeta^{2k}}, &\t{if}\ i=j,\\0&\t{otherwise},\end{cases}$$
$$
	u_{ij}=\frac{(\frac{i}{p})\zeta^{-j-2i}+(\frac{j}{p})\zeta^{-2j-i}}  { \zeta^{-i-j}+(\frac{i}{p})(\frac{j}{p}) }\ \ \t{and}\ \
	v_{ij}=\zeta^{2ij}.$$
Vsemirnov \cite[Theorem 2]{{V1}} proved that
$$ C=\lambda VDUDV,$$	
where $\lambda=(\frac{2}{p})\sqrt{p}\,\zeta^{(p-1)/4}$.

\begin{lemma}\label{New-decompose} Let $p$ be a prime with $p\eq1\pmod4$. For the matrix
$$C(x):=\left[x+\left(\frac{j-i}{p}\right)\right]_{0\le i,j\le \frac{p-1}{2}},
$$
	we have
	\begin{align}\label{cx}
C(x) =\lambda VD\tilde{ U}DV,
	\end{align}
where the top-left entry of $ \tilde{U} $is $ (\frac{2}{p}) \sqrt{p}\, x $ and all the other entries are the same as $U$.
\end{lemma}
\Proof.
We now introduce $n+1$-dimensional vectors
$
\mathbf e_0=(1,0,\ldots,0)^T, 
\mathbf e=(1,\ldots,1)^T.
$ 
Since  $u_{00}=0$,
\[
\widetilde U-U
=
\(\frac{2}{p}\)\sqrt p\,x\,
\mathbf e_0\mathbf e_0^T.
\]
So, 
\[
\lambda VD(\widetilde U-U)DV
=
\lambda
 \(\frac{2}{p}\)
\sqrt p\,x\,V D\mathbf e_0\mathbf e_0^T DV.
\]
Recall that $D$ is diagonal and the first column of
$V $ is  $\mathbf e$. Hence
\[
V D\mathbf e_0
=V(d_{00},0,\ldots,0)^T=d_{00} V \mathbf e_0
=d_{00} \mathbf e.
\]
By symmetry,  $\mathbf e_0^T DV=d_{00}\mathbf e^T$.
Therefore
\[
\lambda VD(\widetilde U-U)DV
=
\lambda
\(\frac{2}{p}\)
\sqrt p\,x\,d_{00}^2
\mathbf e\mathbf e^T.
\]
 Set $n=(p-1)/2$. Then
\begin{align*}
	\frac{1}{	d_{00}^{2} }&= \prod_{0< k\le n }^{}  (1-\zeta^{2k})^2
	 \\
	&=  \prod_{k=1 }^{n}  (-\zeta^{2k})  (1-\zeta^{-2k})  (1-\zeta^{2k})     \\
	&=  (-1)^n\zeta^{n(n+1)}  \prod_{k=1 }^{n}  (1-\zeta^{p-2k}) (1-\zeta^{2k})
\\&=  \zeta^{(p+1)(p-1)/4}   \prod_{r=1 }^{p-1}  (1-\zeta^{r})   \\&= \zeta^{(p-1)/4}\lim_{x\to1}\f{x^p-1}{x-1}=p\zeta^{(p-1)/4}.
\end{align*}
Thus
$d_{00}^{2}=p^{-1} \zeta^{-(p-1)/4}.$
Since $\lambda
=
(\frac{2}{p})
\sqrt p\,
\zeta^{(p-1)/4},$ we obtain that 
$$\lambda
\(\frac{2}{p}\)
\sqrt p\,
d_{00}^2
=
1.$$
So,
\[
\lambda VD(\widetilde U-U)DV
=x \mathbf e\mathbf e^T.
\]
Consequently, 
\[
\lambda VD\widetilde U DV
= \lambda VD U DV +\lambda VD(\widetilde U-U)DV
=C+x \mathbf e\mathbf e^T=C(x).
\]
 This completes the proof.
 \qed

	\section{Proof of \eqref{simple} in the case $p\eq1\pmod4$}
\setcounter{lemma}{0}
\setcounter{theorem}{0}
\setcounter{corollary}{0}
\setcounter{remark}{0}
\setcounter{equation}{0}
\setcounter{conjecture}{0}

Let  $p$ be a prime with $p\equiv 1\pmod4$, and set $n=(p-1)/2$.
Define
	$$R(x):=\left[x+\left(\frac{j-i}{p}\right)\right]_{1\le i,j\le n}.
$$
Since $p\equiv1\pmod4$, we have
$(\frac{-1}{p})=1$, and hence
\[
\left(\frac{j-i}{p}\right)
=
\left(\frac{i-j}{p}\right).
\]
Therefore, $\det R(x)$ is exactly the determinant appearing in
 (\ref{simple}).

Since  $\det R(x)$ is the cofactor of the top-left entry of $C(x)$, we have
$$\det R(x)= ( 1  ,  0,0,\ldots,0  ) \adj (C(x))  ( 1  ,  0,0,\ldots,0  )^T.$$
By Lemma \ref{New-decompose},
\begin{align}
	  \adj (C(x))=\lambda^n\cdot	\adj (V) \cdot	\adj(D) \cdot\adj (\tilde{ U})  \cdot 	\adj (D)\cdot	 \adj (V) . \nonumber
\end{align}
So,
\begin{align}
	 \det R(x)=&\lambda^n\cdot ( 1  ,  0,0,\ldots,0  )\cdot 	\adj (V) \cdot \adj (D) \cdot\adj (\tilde{ U})    \nonumber \\
	 &	\cdot	\adj (D)\cdot	\adj (V)  \cdot  ( 1  ,  0,0,\ldots,0  )^T.
\end{align}

 Recall that $\adj (V)=[V_{ji}]_{0\le i,j\le n}$ with
 $V_{ji}$ the cofactor of $v_{ji}$ in the matrix $V$.
It is easy to see that
\begin{align}
	 ( 1  ,  0,0,\ldots,0  )\cdot 	\adj (V) \cdot	\adj (D)=	\det (D) \( \frac{ V_{00}}{d_{00}} ,   \frac{ V_{10}}{d_{11}}, \frac{ V_{20}}{d_{22}},\ldots, \frac{ V_{n0}}{d_{nn}} \).
\end{align}

Let $0\le k\le n $.  As determinants of Vandermonde's type can be evaluated, we get
	\begin{align}
	V_{k0}&=V_{0k}=(-1)^{0+k}\cdot\det [\zeta^{2ij}]_{1\le i\le n\atop  0\le j\le n, j\ne k}  \nonumber    \\
		&=(-1)^{k}\cdot     \prod_{j=0\atop j\ne k}^{n} \zeta^{2j}  \cdot   \det  [\zeta^{2(i-1)j}]_{1\le i\le n\atop  0\le j\le n, j\ne k}   \nonumber \\
	&= \zeta^{n(n+1)}  \cdot  \zeta^{-2k} \cdot \frac{\prod_{0\le i<j \le n}(\zeta^{2j}-\zeta^{2i})}{\prod_{0\le j\le n\atop j\ne k}(\zeta^{2k}-\zeta^{2j})}. \nonumber
\end{align}
Since
 $$d_{kk} \prod_{0\le j\le n\atop j\ne k}(\zeta^{2k}-\zeta^{2j})=1
 \ \t{and}\ \det V=\prod_{0\le i<j \le n}(\zeta^{2j}-\zeta^{2i}),$$
 we have
 \begin{align}
  \frac{ V_{k0}}{d_{kk}}=   \zeta^{n(n+1)}\cdot \det V\cdot \zeta^{-2k}.\nonumber
 	\end{align}

 Define
 \begin{equation}
 	\boldsymbol{\alpha}:=    ( 1  , \zeta^{-2}, \zeta^{-4},\ldots,\zeta^{-2n} )^T.\nonumber
 \end{equation}
 Then
 $$( 1  ,  0,0,\ldots,0  )\cdot 	\adj (V) \cdot	\adj(D)=	\zeta^{n(n+1)}\cdot \det D\cdot  \det V	 \cdot \boldsymbol{\alpha}^T.$$
 By symmetry,
 $$\adj (D)\cdot	\adj (V) 	( 1  ,  0,0,\ldots,0  )^T=	\zeta^{n(n+1)}\cdot \det D\cdot  \det V\cdot 	 \boldsymbol{\alpha}.$$
Hence
\begin{equation}\label{Rx}
	\det R(x)= p^{n/2}   \zeta^{\frac{n^2}{2 }  +2n(n+1)}\cdot  (\det V)^2 \cdot (\det D)^2\cdot  \boldsymbol{\alpha}^T \cdot   	\adj (\tilde{ U})  \cdot\boldsymbol{\alpha}.
\end{equation}

  Now it remains  to evaluate $\boldsymbol{\alpha}^T \cdot   	\adj (\tilde{ U})  \cdot\boldsymbol{\alpha}$.  By Lemma \ref{mdl}, 
  
  \[
\det (\boldsymbol{\alpha} \boldsymbol{\alpha}^T +\tilde{ U}	)
=
\det\widetilde U
+
\boldsymbol{\alpha}^T\adj(\widetilde U)\boldsymbol{\alpha},
\]
and hence
  \begin{equation}\label{alphaUalpha}
  	\boldsymbol{\alpha}^T \cdot   	\adj (\tilde{ U})\cdot\boldsymbol{\alpha}=\det (\boldsymbol{\alpha} \boldsymbol{\alpha}^T +\tilde{ U}	)- \det \tilde{ U} .
  \end{equation}
Set $x_i=(\frac{i}{p})\zeta^{-i}$ for $1\le i\le n$. Then
   \begin{equation*}
  	\boldsymbol{\alpha}^T=  ( 1  , x_1^{2}, x_2^{2},\ldots,x_n^{2} ).
  \end{equation*}
 Define
\begin{displaymath}
	G =	\left[ \begin{array}{cccccccc}
		1& & \  &  \\
		&	 x_1^{-1}&    &  \\
		& &	 x_2^{-1}&    &  \\
		&	& & \ddots & \\
		&	& &  &	 x_n^{-1}\\
	\end{array} \right],
\end{displaymath}
and write $y=(\frac{2}{p}) \sqrt{p} x $.
Then
\begin{equation}\label{u}
	G  \tilde{ U}	 G =\left[ \begin{array}{c|cccc}
		y & 1& \cdots & 1 \\
		\hline
		1	 &    & \\
		\vdots & &\bigg[\frac{x_i+x_j}{1+x_ix_j}\bigg]_{1\le i,j\le n} \\
		1 &  & \\
	\end{array} \right].
\end{equation}

We next simplify $\det (\boldsymbol{\alpha} \boldsymbol{\alpha}^T +\tilde{ U}	)$.
It is easy to verify that
\begin{align}
&G^2 (\boldsymbol{\alpha} \boldsymbol{\alpha}^T +\tilde{ U}	) G^2  	 \nonumber\\
=	&
\left[ \begin{array}{c|cccc}
		y+1& 1+x_1^{-1}  &  1+x_2^{-1}   & \cdots &  1+x_n^{-1}  \nonumber  \\
		\hline
	 1+x_1^{-1}   & & & \\
	 1+x_2^{-1}  	 & &   & \\
		\vdots & && \bigg[1+\frac{x_i^{-1}+x_j^{-1}}{1+x_ix_j}\bigg]_{1\le i,j\le n} &  \nonumber\\
	 1+x_n^{-1}  & & & \\
	\end{array} \right].
\end{align}
Introduce the bordered matrix obtained by adjoining a new first row and
a new first column, with the first column equal to
$(1,0,\ldots,0)^T$. Expanding its determinant along the first column
shows that its determinant is equal to
$\det\bigl(G^2(\alpha\alpha^T+\widetilde U)G^2\bigr)$.
Subtracting the new first row from each of the remaining rows
does not change the determinant. Thus, we obtain
\begin{align}
	&\ \det (	G^2 (\boldsymbol{\alpha}  \boldsymbol{\alpha}^T +\tilde{ U}	) G^2 )    \nonumber \\
	=&\	\det\left[\begin{array}{cc|cccc}
			1 & 1 &1&1& \cdots & 1  \nonumber \\
		0& y+1  &  1+x_1^{-1}   & 1+x_2^{-1}   & \cdots &  1+x_n^{-1}   \nonumber \\
		\hline
	0& 	1+x_1^{-1}   & & & \\
	0& 	1+x_2^{-1}  	 & &   & \\
\vdots & 	\vdots & & &\bigg[1+\frac{x_i^{-1}+x_j^{-1}}{1+x_ix_j}\bigg]_{1\le i,j\le n}    \nonumber \\
	0& 	1+x_n^{-1}  & & &     \nonumber\\
	\end{array} \right]     \nonumber\\
=&\	\det L,  \nonumber
\end{align}
where
\begin{align}
	 L=\left[  \begin{array}{cc|cccc}
	 	1 & 1 &1& 1&\cdots & 1     \nonumber \\
	 	-1& y  &  x_1^{-1}   &  x_2^{-1}   &\cdots &  x_n^{-1}     \nonumber \\
	 	\hline
	 	-1& 	x_1^{-1}   & &   \nonumber \\
	 	-1& 	x_2^{-1}  	 &    &     \nonumber\\
	 	\vdots & 	\vdots & &&\bigg[\frac{x_i^{-1}+x_j^{-1}}{1+x_ix_j}\bigg]_{1\le i,j\le n}    \nonumber \\
	 	-1& 	x_n^{-1}  &  & \\
	 \end{array} \right].
\end{align}

Define the diagonal matrix
\begin{displaymath}
	\tilde{G }=	\left[ \begin{array}{cccccccc}
		1& & \  &  \\
			&1& & \  &  \\
	&	&	 x_1^{-1}&    &  \\
	&	& &	 x_2^{-1}&    &  \\
	&	&	& & \ddots & \\
	&	&	& &  &	 x_n^{-1}\\
	\end{array} \right].
\end{displaymath}
Then
  \begin{align}
  	\tilde{G }^{-1}L\tilde{G }^{-1}=	\left[\begin{array}{cc|cccc}
  		1 & 1 & x_1&x_2& \cdots & x_n \\
  		-1& y  &1  &1 & \cdots &  1 \\
  		\hline
  		-x_1& 1  & & \\
  		-x_2& 	1	 &    & \\
  		\vdots & 	\vdots & &&\bigg[\frac{x_i+x_j}{1+x_ix_j}\bigg]_{1\le i,j\le n} \\
  		-x_n& 	1 &  & \\
  	\end{array} \right].
  \end{align}
Note that $\det G =\det \tilde{ G}$.
By linearity of the determinant with respect to the first column, writing
\[
\begin{pmatrix}
1\\ -1\\ -x_1\\ \vdots\\ -x_n
\end{pmatrix}
=
\begin{pmatrix}
0\\ -1\\ -x_1\\ \vdots\\ -x_n
\end{pmatrix}
+
\begin{pmatrix}
1\\ 0\\ 0\\ \vdots\\ 0
\end{pmatrix},
\]
we obtain
\begin{align}
	\det (	G (\boldsymbol{\alpha} \boldsymbol{\alpha}^T +\tilde{ U}	) G)\nonumber
=	&\	\det  \tilde{G }^{-1}L\tilde{G }^{-1}	  \nonumber\\
=&\	\det\left[\begin{array}{cc|cccc}
		0 & 1 & x_1&x_2& \cdots & x_n \\
		-1& y  &1  &1 & \cdots &  1\nonumber \\
		\hline
		-x_1& 1  & &\nonumber \\
		-x_2& 	1	 &    & \nonumber \\
		\vdots & 	\vdots & &&\bigg[\frac{x_i+x_j}{1+x_ix_j}\bigg]_{1\le i,j\le n}\nonumber \\
		-x_n& 	1 &  &\nonumber \\
	\end{array} \right]   \nonumber  \\
&\ +\det\left[\begin{array}{cc|cccc}
	1 & 1 & x_1&x_2& \cdots & x_n   \nonumber\\
	0& y  &1  &1 & \cdots &  1    \nonumber\\
	\hline
	0& 1  & &   \nonumber \\
0& 	1	 &    &      \nonumber\\
	\vdots & 	\vdots & &&\bigg[\frac{x_i+x_j}{1+x_ix_j}\bigg]_{1\le i,j\le n}      \nonumber \\
	0& 	1 &  &      \nonumber\\
\end{array} \right]  .      \nonumber\\
\nonumber
\end{align}
Expanding the last determinant along the first column we obtain that 
\begin{equation}\label{alph+u}\begin{aligned}
&\	\det (	G(\boldsymbol{\alpha} \boldsymbol{\alpha}^T +\tilde{ U}	) G) -\det\left[\begin{array}{c|cccc}
	 y & 1& \cdots & 1 \\
	 	\hline
		1	 &    & \\
		\vdots & &\bigg[\frac{x_i+x_j}{1+x_ix_j}\bigg]_{1\le i,j\le n} \\
		1 &  & \\
\end{array} \right]
\\=&\
 -\det\left[\begin{array}{cc|cccc}
	0 & 1 & x_1&x_2& \cdots & x_n        \\
	1& y  &1  &1 & \cdots &  1       \\
	\hline
	x_1& 1  & &      \\
	x_2& 	1	 &    &        \\
	\vdots & 	\vdots & &&\bigg[\frac{x_i+x_j}{1+x_ix_j}\bigg]_{1\le i,j\le n}        \\
	x_n& 	1 &  &        \\
\end{array} \right].         \\
\end{aligned}\end{equation}
Combining (\ref{u}) with (\ref{alph+u}),  we obtain
\begin{equation}\label{GGG}\begin{aligned}
	&\	\det (	G (\boldsymbol{\alpha}  \boldsymbol{\alpha}^T +\tilde{ U}	) G) -\det (G  \tilde{ U}	 G)\\
	=&\	-\det\left[\begin{array}{cc|cccc}
		0 & 1 & x_1&x_2& \cdots & x_n\\
		1& y  &1  &1 & \cdots &  1 \\
		\hline
		x_1& 1  & & \\
		x_2& 	1	 &    & \\
		\vdots & 	\vdots & &&\bigg[\frac{x_i+x_j}{1+x_ix_j}\bigg]_{1\le i,j\le n} \\
		x_n& 	1 &  &
	\end{array} \right].
\end{aligned}\end{equation}
Exchanging the first two rows
and the first two columns does not change the determinant, because the two sign changes cancel. So the last determinant 
\begin{align}
&\ \det\left[\begin{array}{cc|cccc}
	0 & 1 & x_1&x_2& \cdots & x_n \nonumber\\
	1& y  &1  &1 & \cdots &  1\nonumber \\
	\hline
	x_1& 1  & &\nonumber \\
	x_2& 	1	 &    & \nonumber\\
	\vdots & 	\vdots & &&\bigg[\frac{x_i+x_j}{1+x_ix_j}\bigg]_{1\le i,j\le n} \nonumber\\
	x_n& 	1 &  &\nonumber \\
\end{array} \right]	
\\
=&\	\det\left[\begin{array}{cc|cccc}
		y & 1&1  &1 & \cdots &  1 \nonumber \\
		1& 0   & x_1&x_2& \cdots & x_n \nonumber \\
		\hline
	 1  &	x_1& &\nonumber \\
		1	 & 	x_2&    &\nonumber \\
		\vdots & 	\vdots & &&\bigg[\frac{x_i+x_j}{1+x_ix_j}\bigg]_{1\le i,j\le n}\nonumber \\
			1 & x_n&  &\nonumber \\
	\end{array} \right]	
\\
=&\	\det\left[\begin{array}{cc|cccc}
		1 & 1&1  &1 & \cdots &  1 \nonumber \\
		1& 0   & x_1&x_2& \cdots & x_n \nonumber \\
		\hline
		1  &	x_1& &\nonumber \\
		1	 & 	x_2&    &\nonumber \\
		\vdots & 	\vdots && &\bigg[\frac{x_i+x_j}{1+x_ix_j}\bigg]_{1\le i,j\le n}\nonumber \\
		1 & x_n&  & \\
	\end{array} \right]	
\\&\ +	\det\left[\begin{array}{cc|cccc}
		y -1& 1&1  &1 & \cdots &  1 \nonumber \\
		0& 0   & x_1&x_2& \cdots & x_n\nonumber  \\
		\hline
		  &	x_1& & \\
		0	 & 	x_2&    & \\
		\vdots & 	\vdots & &&\bigg[\frac{x_i+x_j}{1+x_ix_j}\bigg]_{1\le i,j\le n} \nonumber\\
		0& x_n&  &\nonumber
	\end{array} \right].	
\end{align}
Expanding the last determinant  along
the first column and using \eqref{GGG}, we obtain
\begin{align*}
	&\	\det (G  \tilde{ U}	 G)-\det (	G  (\boldsymbol{\alpha}  \boldsymbol{\alpha}^T +\tilde{ U}	) G) \\
	=&\ \det\left[\begin{array}{cc|cccc}
		1 & 1&1  &1 & \cdots &  1     \nonumber  \\
		1& 0   & x_1&x_2& \cdots & x_n \nonumber \\
		\hline
		1  &	x_1& &\nonumber \\
		1	 & 	x_2&    &\nonumber \\
		\vdots & 	\vdots & &&\bigg[\frac{x_i+x_j}{1+x_ix_j}\bigg]_{1\le i,j\le n}           \nonumber \\
		1 & x_n&  & \nonumber\\
	\end{array} \right]	 \nonumber
\\&\
	+(y-1) \det \left[ \begin{array}{c|cccc}
		0& x_1 & x_2  & \cdots & x_n\nonumber \\
		\hline
		x_1 & & & \nonumber\\
		x_2	 & &   &\nonumber \\
		\vdots & & &\bigg[\frac{x_i+x_j}{1+x_ix_j}\bigg]_{1\le i,j\le n}\nonumber \\
		x_n& & &\nonumber \\
	\end{array} \right].
\end{align*}

Let $x_{-1}=1$ and $x_{0}=0$. Then
the last two determinants are
 \begin{equation*}
 	\det \bigg[\frac{x_i+x_j}{1+x_ix_j}\bigg]_{-1\le i,j\le n}
 \ \t{and}\
  \det \bigg[\frac{x_i+x_j}{1+x_ix_j}\bigg]_{0\le i,j\le n}
\end{equation*}
respectively.
Thus,
\begin{align}\label{alphaU-U}
	 &\det (G  \tilde{ U}	 G)-\det (	G  (\boldsymbol{\alpha}  \boldsymbol{\alpha}^T +\tilde{ U}	) G)	  \nonumber\\
=	&\det \bigg[\frac{x_i+x_j}{1+x_ix_j}\bigg]_{-1\le i,j\le n} +(y-1)\det \bigg[\frac{x_i+x_j}{1+x_ix_j}\bigg]_{0\le i,j\le n} .
\end{align}
To evaluate the right side of (\ref{alphaU-U}), we define
\begin{align}
	I =\det \bigg[\frac{x_i+x_j}{1+x_ix_j}\bigg]_{-1\le i,j\le n}
\ \t{and}\
	J =\det \bigg[\frac{x_i+x_j}{1+x_ix_j}\bigg]_{0\le i,j\le n} .
\end{align}
Then
\begin{align}\label{I+J}
	\det (	G(\boldsymbol{\alpha}  \boldsymbol{\alpha}^T +\tilde{ U}	)G) -\det (G  \tilde{ U} G)= - I-(y-1)J.
\end{align}
By Lemma \ref{uv} and noting that $x_{-1}=1$ and $x_{0}=0$, we have
\begin{align*}
I=	(-1)^{\frac{p+3}{4}} \prod_{i=1}^{n}x_i^2(  1-x_i )^2
	\cdot\prod_{1\le i<j\le n}^{}  (x_i-x_j)^2 \cdot\prod_{i=1}^{n}  \prod_{j=1}^{n}(1+x_ix_j)^{-1}
\end{align*}
and
\begin{align}
 J	=&\ (-1)^{\f{p-1}4}\cdot\frac{1}{2}\(\prod_{i=1}^{n}(  1+x_i )^2- \prod_{i=1}^{n}(  1-x_i )^2\)   \nonumber\\
	&\ \times\prod_{1\le i<j\le n}^{}(x_i-x_j)^2\cdot\prod_{i=1}^{n}  \prod_{j=1}^{n}(1+x_ix_j)^{-1}\cdot\prod_{i=1}^{n}x_i^2.\nonumber
\end{align}
Set
\begin{equation*}
	f_1= \prod_{1\le i<j\le n}^{}\(   \(\frac{j}{p}\)\zeta^{j}- \(\frac{i}{p}\)\zeta^{i} \)
	= \prod_{1\le i<j\le n}^{}(  x_j^{-1}- x_i^{-1} )
\end{equation*}
and
\begin{equation*}
	f_2= \prod_{1\le i<j\le n}^{}\( 1+  \(\frac{j}{p}\)\zeta^{j}\(\frac{i}{p}\)\zeta^{i} \)	= \prod_{1\le i<j\le n}^{}( 1+ x_j^{-1} x_i^{-1} ).
\end{equation*}
Then one can check that
\[
\frac{f_1^2}{f_2^2}
=
\prod_{1\le i<j\le n}
\frac{(x_i-x_j)^2}{(1+x_ix_j)^2}.
\]
We also use
\[
\prod_{i=1}^{n}\prod_{j=1}^{n}(1+x_ix_j)
=
\prod_{i=1}^{n}(1+x_i^2)
\cdot\prod_{1\le i<j\le n}(1+x_ix_j)^2.
\]
Consequently,
\begin{align}
	\prod_{1\le i<j\le n}^{} (x_i-x_j)^2 \cdot\prod_{i=1}^{n}  \prod_{j=1}^{n}(1+x_ix_j)^{-1}\cdot  \prod_{i=1}^{n}x_i^2= f_1^2\  f_2^{-2} \prod_{i=1}^{n}\l((1+ x_i^{2})^{-1}x_i^2\r).
	\nonumber
\end{align}
By \cite[Corollary 2 (4.6)]{V1}, we have
\begin{align}
	\prod_{i=1}^{n}(1+ x_i^{2})^{-1}\cdot  \prod_{i=1}^{n}x_i^2=\zeta^{-n(n+1)/2} \(\frac{2}{p}\).  \nonumber
\end{align}
So,
\begin{align*}
	\prod_{1\le i<j\le n}^{} (x_i-x_j)^2\cdot\prod_{i=1}^{n}  \prod_{j=1}^{n}(1+x_ix_j)^{-1}\cdot  \prod_{i=1}^{n}x_i^2= f_1^2 f_2^{-2} \zeta^{-n(n+1)/2} \(\frac{2}{p}\).	\nonumber\\
\end{align*}
Then, in light of Lemma \ref{a+-b}, we have
$$
I	=  ( a_p' \sqrt{p}-b_p' p)  f_1^2  f_2^{-2} \zeta^{-n(n+1)} \(\frac{2}{p}\)$$
and
$$
	J	=  a_p' \sqrt{p}\, f_1^2  f_2^{-2} \zeta^{-n(n+1)} \(\frac{2}{p}\).$$
 Combining the above results with (\ref{I+J}),  and noting that $y= (\frac{2}{p})\sqrt{p}\,x$
 and  $	(\det G)^2= \zeta^{n(n+1)} $, we obtain
\begin{align*}
&\ \det (\boldsymbol{\alpha}  \boldsymbol{\alpha}^T +\tilde{ U}	)- \det \tilde{ U}
\\=&\  \zeta^{-n(n+1)}(-I-(y-1)J)\\
=&\  \(\(\frac{2}{p}\)b_p'-a_p'x \) p\zeta^{-2n(n+1)}   f_1^2  f_2^{-2}.
\end{align*}
This, together with  (\ref{Rx})  and (\ref{alphaUalpha}), yields that
\begin{align*}
	 \det R(x)
	 =&\ \(\(\frac{2}{p}\)b_p'-a_p'x \) \zeta^{\frac{(p-1)^2}{8}}  p^{\frac{p+3}{4}}
 (\det V)^2  (\det D)^2 f_1^2  f_2^{-2}      \\
	 =&\ \l(\frac{2}{p}\r)b_p'-a_p'x \quad \ (\t{by \cite[(4.10)]{V1}}).
\end{align*}
This completes our proof of \eqref{simple} for $p\eq1\pmod4$. \qed

\medskip
\noindent{\bf Statements and Declarations}. There are no competing interests. This original paper contains no data, and  it has not been submitted elsewhere.

\noindent{\bf Acknowledgement}. We would like to thank the referee for the helpful comments.

\end{document}